\documentclass[a4paper,11pt,preprint]{article}

\usepackage{geometry}
\usepackage{amsmath}
\usepackage{amsfonts}
\usepackage{amsthm}
\usepackage{amssymb}
\usepackage{upref}
\usepackage{color}
\usepackage{graphicx}
\usepackage{hyperref}
\usepackage{array}
\usepackage[toc,page]{appendix}
\usepackage{mathtools,tikz-cd}

\graphicspath{{../figures/}{../figures_eps/}}

\theoremstyle{plain}
\newtheorem{thm}{Theorem}[section]
\newtheorem{lemma}[thm]{Lemma}
\newtheorem{cor}[thm]{Corollary}
\newtheorem{prop}[thm]{Proposition}

\newtheorem{rem}[thm]{Remark}
\theoremstyle{definition}
\newcommand{\Per}{\operatorname{Per}}
\newcommand{\ds}{\displaystyle}
\newcommand{\di}{\operatorname{div}}

\usepackage{indentfirst} 
\newenvironment{ack}{{\bf Acknowledgements.}}

\newcommand{\wconv}{\rightharpoonup}

\date{}

\begin{document}

\title{Regularity result for a shape optimization problem under perimeter constraint}
\author{Beniamin \textsc{Bogosel}}

\newcommand{\Addresses}{{
  \bigskip
  \footnotesize

  B.~Bogosel, \textsc{Laboratoire Jean Kuntzmann, Grenoble, France}\par\nopagebreak
  \textit{E-mail address}: \texttt{beniamin.bogosel@cmap.polytechnique.fr}

}}


\maketitle

\begin{abstract}
We study the problem of optimizing the eigenvalues of the Dirichlet Laplace operator under perimeter constraint. We prove that optimal sets are analytic outside a closed singular set of dimension at most $d-8$ by writing a general optimality condition in the case the optimal eigenvalue is multiple. As a consequence we find that the optimal $k$-th eigenvalue is strictly smaller than the optimal $(k+1)$-th eigenvalue. We also provide an elliptic regularity result for sets with positive and bounded weak curvature. 
\end{abstract}



\section{Introduction}
For a given domain $\Omega \in \Bbb{R}^d$ we can consider the eigenvalue problem for the Laplace operator with Dirichlet boundary conditions defined by
\[ \begin{cases}
-\Delta u = \lambda_k(\Omega) u & \text{ in } \Omega \\
u = 0 & \text{ on }\partial \Omega.
\end{cases} \] This operator has compact resolvent if $L^2(\Omega)$ injects compactly in $H^1(\Omega)$ and in this case the spectrum consists of an increasing sequence of eigenvalues
\[ 0<\lambda_1(\Omega) \leq \lambda_2(\Omega) \leq ... \leq \lambda_k(\Omega) \leq ... \to \infty.\] 
The optimization of the eigenvalues of the Dirichlet Laplace operator under various constraints is a classical problem which recently had important developments. Authors were initially concerned with considering volume constraints, but in some recent works like \cite{bucbuthenper}, \cite{deveper} the perimeter constraint gained some interest. In fact, having a perimeter constraint allows the use of techniques related to perimeter quasi-minimizers \cite{tamanini_reg} which allow important gains in regularity properties. More precisely, the basic form of the problem which we consider is the minimization of $\lambda_k(\Omega)$ under a perimeter constraint
\[ \min_{\Per(\Omega) = c} \lambda_k(\Omega).\]
Using the nice behavior of these eigenvalues under homotheties, $\lambda_k(t\Omega) = \lambda_k(\Omega)/t^2$ one can deduce that the above problem is equivalent, up to a homothety, to the problem
\begin{equation} \min_{\Omega \in \Bbb{R}^d} \left( \lambda_k(\Omega)+\Per(\Omega)\right).
\label{persum}
\end{equation}
This latter formulation has the advantage of being unconstrained and will be used throughout the paper. In \cite{ldpv} the authors study existence and regularity properties of problems of the form 
\[ \min_{|\Omega|=m}\left( \Per(\Omega)+\mathcal{G}(\Omega)\right)\] where $\mathcal{G}$ is the Dirichlet energy or a spectral functional. The authors focus on showing how to recover the $C^{1,\alpha}$ regularity of the relaxed optimizers. They also describe how to perform a bootstrap argument in certain cases and be able to show that optimal sets are $C^\infty$. The case where the functional depends on spectral quantities is not treated in full detail as one needs to take care of the case where the optimal eigenvalues are multiple, thus losing some differentiability properties.

The aim of this paper is to improve the known regularity properties for problem \eqref{persum}. In fact for the same problem with volume constraint regularity issues are more delicate. Existence results are proved in \cite{bucur-mink}, \cite{maz_prat} in the class of quasi-open sets, which are basically level sets of $H^1$ functions. Initial regularity results are presented in \cite{bmpv} where the authors prove that in certain cases the solution to the analogue spectral optimization problem under volume constraint is an open set. From now on we consider only the perimeter constraint.

The Faber-Krahn inequality and the isoperimetric inequality immediately imply that the first eigenvalue is minimized by the ball in any dimension. A first result concerning higher eigenvalues is due to Bucur, Buttazzo and Henrot \cite{bucbuthenper}. The authors consider the minimization of the second eigenvalue among domains $\Omega$ in $\Bbb{R}^2$. Some results particular to this case, like convexity of the optimizers, $H^2$ regularity for elliptic problems on convex domains \cite{grisvard} and simplicity of the second eigenvalue, allow the use of a bootstrap procedure which shows that the optimal set is $C^\infty$.

In \cite{deveper} De Philippis and Velichkov show that problem \eqref{persum} has solutions in every dimension and the optimal shapes have regularity at least $C^{1,\alpha}$ outside a closed set of Hausdorff dimension less or equal than $d-8$. The fact that minimizers are not known and do not have a simple structure even in the case of the second eigenvalue \cite{bucbuthenper} motivated the numerical study of problem \eqref{persum}. Simulations were performed in $\Bbb{R}^2$ and $\Bbb{R}^3$ with two distinct methods in \cite{antunesf-per} and \cite{boouper}. In addition to numerical simulations Bogosel and Oudet \cite{boouper} also present an optimality condition applicable in the case of multiple eigenvalues inspired by methods from \cite{iliasoptim}. The need of these more general optimality conditions is motivated by the numerical observations. Indeed, we observe that except some particular cases the eigenvalues of the optimal shapes seem to have multiplicity higher than one. Thus classical methods based on optimality conditions obtained directly from the derivative of the eigenvalue, like the one used in \cite{bucbuthenper}, are not always applicable. The more general optimality condition presented in \cite{boouper} states that if the optimal set $\Omega^*$ is of class $C^3$ then there exist a family of eigenfunctions $\{\phi_i\}_{i=1}^m$ such that the mean curvature $\mathcal{H}$ of  $\Omega^*$  can be expressed as 
\[ \mathcal{H} =\sum_{i=1}^m (\partial_n \phi_i)^2.\]
This optimality condition allows the use of a bootstrap argument which shows that optimal shapes are smooth. These optimality conditions were also used in \cite{boouper} as a tool for validating the numerical computations.

The purpose of this note is to remove the additional $C^3$ regularity hypothesis used in the deduction of the above optimality condition. We prove the following result.

{\bf Theorem \ref{main-regularity}}\emph{ If $\Omega^*$ is a local minimizer for problem \eqref{persum} then $\Omega^*$ is analytic for $d\leq 7$. If $d\geq 8$ then $\Omega^*$ is analytic outside a closed set of Hausdorff dimension $d-8$.}

The arguments used in the proof are different from the ones used in \cite{buboou} which used heavily the $C^3$ regularity. There are also common points such as the analytic parametrization of eigenvalues under regular perturbations and the use of a variant of the Hahn-Banach theorem to conclude. A summary of the steps employed in the proof is presented below.
\begin{enumerate}
\item We want to be able to recover some differentiability information in the case the eigenvalue is multiple. It is well known (see \cite{henrot-pierre}) that in this case the eigenvalue may not be differentiable. Nevertheless, applying a classical result of Rellich \cite{rellich-perturbation} we can deduce that the eigenvalues and eigenfunctions of an analytic perturbation of the domain can be reparametrized analytically. This fact implies, in particular, that even if we lose differentiability we still recover some facts about directional derivatives.
\item A second aspect is writing the lateral derivatives of the eigenvalues using expressions similar to the shape derivative formula. In order to do this we would like to have more information about the regularity of the eigenfunctions. In the work of Grisvard \cite[Theorem 3.2.1.2]{grisvard} it is proved that if $\Omega$ is convex then solving $-\Delta u = f$ in $H_0^1(\Omega)$ with $f \in L^2(\Omega)$ implies that $u$ is a $H^2(\Omega)$ function. The proof of this uses the fact that convex sets can be approximated uniformly from inside with $C^2$ convex sets. This brings us naturally to the next point.

\item In \cite{deveper} it is proved that solutions of problem \eqref{persum} exist, are of class $C^{1,\alpha}$ and that any optimizer has non-negative and bounded weak mean curvature. The question which arises is the following: can we approximate $C^{1,\alpha}$ sets with weak non-negative bounded mean curvature with smoother sets (at least $C^2$) with non-negative mean curvature? The answer is affirmative and a proof is presented in \cite[Lemma 3.8]{metzger_schulze}

\item The final step consists in using a variant of the Hahn-Banach theorem which allows us to deduce that the mean curvature is a convex combination of squares of normal derivatives of an orthonormal basis of the eigenspace corresponding to the $k$-th eigenvalue. Once we have this formula we use a result of Landais \cite{landais} on the regularity of solutions to the mean curvature equation to deduce that the optimal set is $C^{2,\alpha}$. Results concerning the regularity of solutions to equations of the mean curvature type can also be found in \cite[Section 7.7]{ambrosiofuscopallara}. A  bootstrap algorithm similar to the one used in \cite{bucbuthenper} or \cite{landais} allows us to conclude that optimal shapes are $C^\infty$. \end{enumerate}

After proving the regularity result we turn our attention to the case where the eigenvalues of the optimal set are multiple. We prove a result which was observed in the numerical computations from \cite{boouper}. The result says that when the optimal $k$-th eigenvalue is multiple then the multiplicity cluster ends at $\lambda_k$, namely $...=\lambda_{k-1}(\Omega^*) = \lambda_k(\Omega^*)<\lambda_{k+1}(\Omega^*)$. This implies, in particular, that the optimal costs 
\begin{equation} c_k = \min_{\Omega \in \Bbb{R}^n} \left(\lambda_k(\Omega)+\Per(\Omega)\right), \ \ d_k = \min_{\Per(\Omega) = 1} \lambda_k(\Omega)
\label{costk}
\end{equation}
form a strictly increasing sequence.

\section{Elliptic regularity for sets of positive mean curvature}

As the results of \cite{deveper} prove that the solutions of problem \eqref{persum}, denoted by $\Omega^*$ in the sequel, are $C^{1,\alpha}$ domains (outside a residual closed set of dimension at most $d-8$) with non-negative weak mean curvature, we start from this setting and we prove the $H^2$ regularity of the eigenfunctions for the optimizers of \eqref{persum}. In the first part of this section we restrict ourselves to dimension $d$ at most $7$, where we do not have singularities. In this case we conclude that $u \in H^2(\Omega)$. At the end of the section we give an argument which allows us to deduce local $H^2$ regularity outside the singular set. This type of results are new up to the author's knowledge. Results relating the regularity of the solution to the curvature of the boundary can be found, for example in \cite[Section 14.3]{gilbarg-trudinger}, however the problem we are interested in is not covered by these. The results in \cite[Chapter 3]{grisvard} concerning elliptic problem on convex domains are similar to the ones we are interested in and we follow similar methods in the proofs below.

 Thus, following the results of \cite{deveper} we know that solutions of \eqref{persum} can be represented locally around each point in the regular part of the boundary as epigraphs of $C^{1,\alpha}$ functions. Thus, in a local coordinate system with origin in a regular point $x_0$, the boundary $\partial \Omega^*$ can be represented as the graph of $h: D \subset \Bbb{R}^{d-1} \in \Bbb{R}$, where $D$ is an open set containing $0$. Moreover, following the results of \cite{deveper}, for every $\psi \in C^1(D),\ \psi\geq 0$ we have
\[ 0\leq  \int_D \frac{\nabla h \cdot \nabla \psi}{\sqrt{1+|\nabla h|^2}} \leq K_\infty  \]
where $K_\infty$ is a constant. Thus the mean curvature of an optimizer $\Omega^*$ can be represeted as an $L^2(\partial \Omega)$ integrable function, since $\Omega^*$ has finite perimeter. Note that if $\Omega^*$ is $C^2$ then these inequalities simply say that its mean curvature is non-negative and bounded above by $K_\infty$. 

In \cite[Chapter 3]{grisvard} we can find results concerning the regularity of elliptic problems in terms of the mean curvature of the boundary. In particular, for convex sets the solution of $-\Delta u = f$ in $H_0^1$ for $f \in L^2$ is an $H^2$ function. The scope of this section is to extend such results to the case of $C^{1,\alpha}$ domains which have weak non-negative mean curvature. In order to do this we follow a similar approach. In a first stage we approximate $C^{1,\alpha}$ domains with non-negative weak mean curvature with $C^2$ sets with non-negative mean curvature. This approximation stage is presented in \cite{metzger_schulze} and uses notions related to  weak mean curvature flows. A second stage, which is a simple consequence of Theorem 3.1.2.1 and Remark 3.1.2.2 in \cite{grisvard} gives some results for $C^2$ domains with non-negative mean curvature. We conclude the general case using an approach similar to Theorem 3.1.2.1 from \cite{grisvard}.

\begin{prop}
Let $\Omega$ be a $C^{1,\alpha}$ domain with non-negative weak mean curvature in $L^2(\partial \Omega)$. Then for each $\varepsilon>0$ there exists a smooth set $\Omega_\varepsilon$, with non-negative mean curvature such that $\Omega_\varepsilon \subset \Omega$ and $d_H(\partial \Omega,\partial \Omega_\varepsilon)<\varepsilon$ (where $d_H$ represents the Hausdorff distance).
\label{pc-approx}
\end{prop}

\emph{Proof:} We see that the set $\Omega$ verifies the hypotheses of Lemma 3.8 in \cite{metzger_schulze}. Indeed, $\Omega$ has non-negative weak mean curvature in $L^2(\partial \Omega)$ so there exists a sequence $(\Omega_\varepsilon)_{\varepsilon>0}$ such that $\Omega_\varepsilon \to \Omega$ in $C^1$ (and thus uniformly), $\Omega_\varepsilon$ are smooth for $\varepsilon>0$ and the mean curvature of $\Omega_\varepsilon$ is strictly positive. \hfill $\square$

We state below a few standard regularity results for regular domains. Proofs can be found, for example, in Theorems 2.2.2.3 and 3.1.2.1 from \cite{grisvard}. \hfill $\square$
\begin{prop}
1. If $\Omega$ is a domain with $C^2$ boundary. Then for every $f \in L^2(\Omega)$ the problem
\[\begin{cases}
-\Delta u = f & \text{ in }\Omega \\
u = 0 & \text{ on }\partial \Omega
\end{cases}\] has a unique solution solution $u \in H^2(\Omega)$.

2. If $\Omega$ is an open, bounded subset of $\Bbb{R}^2$ with $C^2$ boundary and non-negative mean curvature, then there exists a constant $C(\Omega)$, which depends only on the diameter of $\Omega$, such that
\[\|u\|_{H^2(\Omega)} \leq C(\Omega) \|\Delta u\|_{L^2(\Omega)}, \]
for every $u \in H^2(\Omega) \cap H_0^1(\Omega)$.
\label{h2-reg}
\end{prop}
Now we are ready to prove the desired regularity result.

\begin{thm}
Suppose $\Omega$ is a $C^{1,\alpha}$ set with non-negative weak mean curvature in $L^2(\partial \Omega)$. Then for every $f \in L^2(\Omega)$ the problem \[\begin{cases}
-\Delta u = f & \text{ in }\Omega \\
u = 0 & \text{ on }\partial \Omega
\end{cases}\] has a unique solution $u \in H^2(\Omega)$.
\label{h2-nreg}
\end{thm}

\emph{Proof:} Proposition \ref{pc-approx} implies the existence of a sequence of $C^2$ sets $\Omega_n$ such that $\Omega_n \subset \Omega$, $\Omega_n$ has non-negative mean curvature and $d_H(\partial \Omega_n,\partial \Omega) \to 0$. We consider the solution $u_m \in H^2(\Omega_m)$ of the Dirichlet problem in $\Omega_m$:
\[\begin{cases}
-\Delta u_m = f & \text{ in }\Omega_m \\
u_m = 0 & \text{ on }\partial \Omega_m,
\end{cases}\]
solution which exists by Proposition \ref{h2-reg}. Applying Theorem 1.5.1.5 from \cite{grisvard} we deduce that $\tilde u_m \in H^1(\Bbb{R}^d)$, where $\tilde u_m$ is the extension by zero of $u_m$ outside $\Omega_m$. Proposition \ref{h2-reg} implies that there exists a constant $C$ such that $\|u_m\|_{H^2(\Omega_m)} \leq C$. This implies that $\tilde u_m$ is bounded in $H^1(\Bbb{R}^d)$ and that the sequences $v_{m,i,j} = \partial_i \partial_j u_m$ are bounded in $L^2(\Bbb{R}^d)$. Therefore, up to choosing a subsequence, we can assume that there exist $U \in H^1(\Bbb{R}^d)$ and $V_{i,j} \in L^2(\Bbb{R}^d)$ such that
\[ u_m \wconv U \text{ weakly in } H^1(\Bbb{R}^d) \text{ and } v_{m,i,j} \wconv V_{i,j} \text{ weakly in }L^2(\Bbb{R}^d).\]

In the following we denote by $u$ the restriction of $U$ to $\Omega$. Since the supports of $\tilde u_m$ are all contained in $\overline \Omega$, it follows that $U$ is also supported on $\overline \Omega$. Thus $u = 0$ on $\partial \Omega$. Now let $\varphi \in C_c^\infty(\Omega)$ be a smooth function with compact support in $\Omega$. Then there exists $m_0$ such that for every $m \geq m_0$ the support of $\varphi$ is contained in $\Omega_m$. Thus, for all $m \geq k_0$ we have
\[ \int_\Omega f\varphi = -\int_{\Omega_m} \Delta u_m \varphi =  \int_{\Omega_m} \nabla u_m \cdot \nabla \varphi = \int_\Omega \nabla \tilde u_m \cdot \nabla \varphi. \]
Taking the limit as $m \to \infty$ we get
\[ \int_\Omega f\varphi = \int_\Omega \nabla u\cdot \nabla \varphi.\]
Since this is true for every $\varphi \in C_c^\infty(\Omega)$ we conclude that $-\Delta u = f$ on $\Omega$ in the sense of distributions. Until now we have the existence of a solution $u \in H^1(\Omega)$. In order to complete the proof we need to prove that the second derivatives of $u$ are in $L^2$. We take again $\varphi \in C_c^\infty(\Omega)$ and for $m \geq m_0$ we have
\[ \int_\Omega \tilde u_m \partial_i\partial_j \varphi = \int_{\Omega_m} u_m \partial_i\partial_j \varphi = \int_{\Omega_m} \partial_i \partial_j u_m \varphi =\int_\Omega v_{m,i,j} \varphi. \]
Taking the limit as $m \to \infty$ we obtain
\[ \int_\Omega u\partial_i\partial_j \varphi =\int_\Omega V_{i,j} \varphi,\]
which means that the distributional derivatives of order $2$ denoted by $\partial_i\partial_j u$ are given by $V_{i,j}$ and they are in $L^2$ for all $i,j = 1...d$. Thus we have proved the existence of a solution $u \in H^2(\Omega)\cap H_0^1(\Omega)$. The uniqueness follows at once using the estimates in Proposition \ref{h2-reg} along with an approximation argument. \hfill $\square$

We note that the above arguments work if $\Omega$ does not have a singular part $\Sigma$. Indeed, the result stated in Proposition \ref{pc-approx} is for the case when $\Omega$ is globally $C^1$. When we have a singular set $\Sigma$ it is no longer possible to use this. We can however prove that outside the singularity the eigenfunction is locally $H^2$. In order to do this, let's recall the following integration by parts formula. 
 \[ \int_\Omega \varphi |D^2u|^2dx + \int_{\partial \Omega} \varphi \mathcal{H}|\nabla u|^2d\sigma = \int_\Omega \varphi (\Delta u)^2dx +\int_\Omega(\nabla\varphi \cdot \nabla u)\Delta u dx-\int_\Omega \nabla \varphi \cdot D^2u\nabla u dx. \]
For the sake of completness an idea of the proof, as well as some references are given in the Appendix \ref{int-parts}. Note that for a solution of \eqref{persum} we know that the corresponding eigenfunction $u$ is Lipschitz so $\nabla u \in L^\infty(\Omega)$. See \cite{deveper} for details. Moreover, since $-\Delta u = \lambda_k u$ and $u \in L^\infty(\Omega)$ we may as well assume that $\Delta u \in L^\infty(\Omega)$. See for example \cite[Example 2.1.8]{davies} where we have the estimate $\|u\|_\infty \leq e^{1/8\pi} \lambda_k(\Omega)^{d/4}$. 
 
Let $x_0$ be a point outside the singular set $\Sigma$. Since this set is closed, there is an open neighborhood of $x_0$ which does not intersect the singular part. Let $B_r$ be a ball centered in $x_0$ of radius $r$ small enough such that $B_r \cap \partial \Omega$ is $C^{1,\alpha}$. Consider a smooth cutoff function $0\leq \varphi \leq 1$ such that, $\varphi = 1$ on $B_{r/2}$ and $\varphi = 0$ on $\partial B_r$, where $B_{r/2}$ is a ball concentric with $B_r$, having radius $r/2$. With these considerations note that 
\[ \int_{\partial (\Omega \cap B_r)} \varphi \mathcal{H}|\nabla u|^2d\sigma \geq 0\]
since on $\partial \Omega \cap B_r$ the curvature is non-negative in the distributional sense and on $\partial B_r \cap \Omega$ we have $\varphi=0$. Thus, applying the above integration by parts formula for $\Omega \cap B_r$ and $\varphi$ described above, we get
\[ \int_{\Omega \cap B_r} \varphi |D^2u|^2dx \leq  \int_{\Omega \cap B_r} \varphi (\Delta u)^2dx +\int_{\Omega\cap B_r}(\nabla\varphi \cdot \nabla u)\Delta u dx-\int_{\Omega \cap B_r} \nabla \varphi \cdot D^2u\nabla u dx. \]
Applying the inequality $ab\leq \frac{a^2+b^2}{2}$ and Cauchy-Schwarz we obtain
 \begin{multline*} \int_{\Omega\cap B_r} \varphi |D^2u|^2 dx \leq \int_{\Omega\cap B_r} \varphi(\Delta u)^2dx+\frac{1}{2}\int_{\Omega \cap B_r} (\Delta u)^2 dx +\frac{1}{2}\int_{\Omega\cap B_r} (\nabla \varphi \cdot \nabla u)^2dx\\
 +{ \int_{\Omega\cap B_r} \nabla \varphi\cdot  D^2u\nabla udx}. 
 \label{ugly}
  \end{multline*}
 Integrating again by parts the last term in the above inequality, we obtain
 \begin{align*}
  \int_\Omega \nabla \varphi\cdot  D^2u\nabla udx =\sum_{i,j=1}^n \int_\Omega \partial_i \varphi \partial_i\partial_j u \partial_j u  \stackrel{\partial_j} = & \\
   \sum_{i,j=1}^n\int_{\partial \Omega} \partial_i u \partial_j u \partial_i \varphi \nu_j d\sigma - \sum_{i,j=1}^n \int_\Omega \partial_i u \partial_i \partial_j \varphi \partial_j u dx - \sum_{i,j=1}^n\int_\Omega \partial_i u \partial_i \varphi \partial_j^2 u dx = \\
    \int_{\partial \Omega} (\nabla u \cdot \nabla \varphi)( \nabla u \cdot n) d\sigma - \int_\Omega \nabla u \cdot D^2\varphi \nabla u dx - \int_\Omega ( \nabla u \cdot \nabla \varphi )\Delta u dx.&
 \end{align*}
Since $\nabla u$ is in $L^\infty(\Omega)$ and in our case $-\Delta u=\lambda_k u \in L^\infty(\Omega)$, the above expression is bounded and we may conclude that 
\[ \int_{\Omega \cap B_r} \varphi |D^2 u|^2dx <+\infty\]
which means that $u \in H^2(B_{r/2})$. Since this argument is valid for any point $x_0$ outside $\Sigma$, it follows that $u \in H^2_{\text{loc}}(\Omega \setminus \Sigma)$.

\section{Analytic perturbations for eigenvalue problems}
\label{perturbations}

We know that the eigenvalues of the Dirichlet Laplace operator on $\Omega$ are differentiable if  they are simple. We refer to \cite{henrot-pierre} for further details. However, numerical results found in \cite{antunesf-per,boouper} show that optimizers of problem \eqref{persum} tend to have multiple eigenvalues. 
It is possible to recover some information regarding the differentiability of the eigenvalues, even when the multiplicity is greater than one, by using the theory of operator perturbations.


Rellich proved in \cite{rellich-perturbation} that analytic perturbations of a self-adjoint operator with the same domain allows an analytic parametrization of the eigenvalues and the eigenfunctions. More recent results can be found  in \cite{perturbation-article}. The article of Micheletti \cite{micheletti} contains details about the application of Rellich's result to the study of eigenvalue problems. Indeed, it is proved that the Laplace operator defined on $\Omega_\varepsilon=(\text{Id}+\varepsilon V)(\Omega)$ has the same spectrum as a self-adjoint operator $A_\varepsilon$ which depends analytically on $\varepsilon$ and has a domain of definition independent of $\varepsilon$.
The result of 
Rellich, stated above, implies that the eigenvalues and eigenvectors of $(-\Delta)$ on $\Omega_\varepsilon$ can be parametrized analytically with respect to $\varepsilon$. Thus, if $\lambda$ is an eigenvalue of $\Omega$ with multiplicity $m$, with the associated orthonormal basis $(\mathfrak{u}_i)_{i=1}^m$, there exist $m$ functions $\varepsilon \mapsto \lambda_{\varepsilon,i}$ and the associated family of eigenfunctions $\varepsilon \mapsto \varphi_{\varepsilon,i}$ orthonormal in $L^2(\Omega_\varepsilon)$, such that both dependencies are analytic in $\varepsilon$, $\lambda_{0,i} = \lambda$, $\varphi_{0,i} = \mathfrak{u}_i$ and
\[ \begin{cases}
 -\Delta \varphi_{\varepsilon,i} = \lambda_{\varepsilon,i} \varphi_{\varepsilon,i} & \text{ in }\Omega_\varepsilon \\
 \varphi_{\varepsilon,i} = 0 & \text{ on }\partial \Omega_\varepsilon,
\end{cases}\]
for every $i=1,...,m$. Differentiating with respect to $\varepsilon$ we obtain the classical derivative formula
\begin{equation}
\frac{d\lambda_{\varepsilon,i}}{d\varepsilon}\bigg|_{\varepsilon = 0}= -\int_{\partial \Omega} \left(\frac{\partial \mathfrak{u_i}}{\partial n}\right)^2V.n d\sigma 
\label{derivatives}
\end{equation}
For details see \cite{oudeteigs} where a formal argument shows how to deduce relation \eqref{derivatives}. The formal argument is rigorous once we know that the eigenvalues and eigenfunctions depend analytically on $\varepsilon$. The $H^2$ regularity (local $H^2$ when singularities are present) of the eigenfunctions proved in the previous section allows writing the derivative of the eigenvalues using the above boundary integral depending on the normal derivative. For details see \cite[Theorem 5.7.4]{henrot-pierre}.

\section{Main result}

\begin{thm}
If $\Omega^*$ is a local minimizer for problem \eqref{persum} then $\Omega^*$ is analytic for $d\leq 7$. If $d\geq 8$ then $\Omega^*$ is analytic outside a closed set of Hausdorff dimension at most $d-8$.
\label{main-regularity}
\end{thm}

\emph{Proof:} Since $\Omega^*$ has $C^{1,\alpha}$ regularity outside an eventual closed singular set $\Sigma$ of dimension at most $d-8$, it is possible to represent the boundary around a point $x_0\in \partial \Omega^*\setminus \Sigma$ as the graph of a $C^{1,\alpha}$ function $h : B_a \subset \Bbb{R}^{d-1} \to \Bbb{R}$, where $B_a$ is a ball of radius $a$, small enough such that we are still outside the singular set. We denote by $\Omega_\varepsilon=(\text{Id}+\varepsilon V)(\Omega^*)$ for a $C^1(\Bbb{R}^d,\Bbb{R}^d)$ vector field. Using this type of graph representation we deduce that if $\psi$ is the perturbation of this graph associated to $V$ then $\varepsilon \mapsto \Per(\Omega_\varepsilon)$ is differentiable at zero and its derivative is equal to 
\[ \int_{B_a} \frac{\nabla \psi \cdot \nabla h}{\sqrt{1+|\nabla h|^2}}.\]
 As noted in the previous section, the function $\varepsilon \mapsto \lambda_k(\Omega_\varepsilon)$ has left and right derivatives. Furthermore, we know that there is an orthonormal basis of the associate eigenspace, denoted $(\mathfrak u_i)_{i=1}^p$, such that the directional derivatives of $\lambda_k$ are in the set
\[ \left( -\int_{B_a} |\nabla \mathfrak{u}_i|^2 \psi(x)dx \right)_{i=1}^p.\]
This is just a simple consequence of the derivative formula \eqref{derivatives} in the local coordinate system. Indeed, the normal derivative is just the gradient since the tangential derivative is zero. On the other hand the term $V.n$ is equal to $\psi(x)/\sqrt{1+|\nabla h|^2}$ and the Jacobian is $\sqrt{1+|\nabla h|^2}$.

 Thus, at the optimum, regardless of the multiplicity and the $C^1$ perturbation considered, there exist left and right derivatives for $\lambda_k(\Omega_\varepsilon)+\Per(\Omega_\varepsilon)$ at $\varepsilon=0$. The case where the eigenvalue is simple is straightforward and is similar to the approach used in \cite{bucbuthenper}. Suppose now that $\lambda_k$ is multiple. The local optimality of $\Omega^*$ implies that the left and right derivatives of $\lambda_k(\Omega_\varepsilon)+\Per(\Omega_\varepsilon)$ at $\varepsilon=0$ in any direction given by $\psi$ in the local coordinates are of different signs. As a consequence, given a perturbation $\psi \in C^1(B_a)$ we have two indices $i,j$ such that 
\begin{equation} \int_{B_a} \frac{\nabla \psi \cdot \nabla h}{\sqrt{1+|\nabla h|^2}}-\int_{B_a} |\nabla \mathfrak{u}_i|^2 \psi \geq 0 
\label{upper}
\end{equation}
\begin{equation} \int_{B_a} \frac{\nabla \psi \cdot \nabla h}{\sqrt{1+|\nabla h|^2}}-\int_{B_a} |\nabla \mathfrak{u}_j|^2 \psi \leq 0 
\label{lower}
\end{equation}
We define the linear functionals $\mathcal{F},\mathcal{G}_i$ on $C^1(B_a)$ by
\[ \mathcal{F}(\psi) = \int_{B_a} \frac{\nabla \psi \cdot \nabla h}{\sqrt{1+|\nabla h|^2}}\]
\[ \mathcal{G}_i(\psi) = \int_{B_a} |\nabla \mathfrak{u}_i|^2 \psi \]
Equations \eqref{upper} and \eqref{lower} tell us that for every $\psi \in C^1(B_a)$ there exist indices $i,j$ such that $\mathcal{F}(\psi) \geq \mathcal{G}_i(\psi)$ and $\mathcal{F}(\psi) \leq \mathcal{G}_j(\psi)$. This implies, in particular, that for every $\psi\in C^1(B_a)$ we have $\mathcal{F}(\psi) \in \text{conv}\{ \mathcal{G}_i(\psi)\}$. Using a variant of the Hahn-Banach separation theorem presented in Proposition \ref{hb-variant} it follows that $\mathcal{F} \in \text{conv} \{ \mathcal{G}_i\}$. For the sake of completeness a proof of the result can be found in Appendix \ref{hb}. For more details see \cite[Section 2.4]{clarke}.  

Therefore there exist $\mu_1,...,\mu_m \in [0,1]$ with $\mu_1+...+\mu_n = 1$ such that
\[ \int_{B_a} \frac{\nabla \psi \cdot \nabla h}{\sqrt{1+|\nabla h|^2}} = \sum_{i=1}^m \mu_i \int_{B_a} |\nabla \mathfrak{u}_i|^2 \psi,\]
for every $\psi \in C^1(B_a)$. Since $h$ is already $C^{1,\alpha}$, we know that each $u_i$ is also $C^{1,\alpha}$ \cite[Theorem  6.18]{gilbarg-trudinger}. Therefore $\nabla \mathfrak u_i$ are all in $C^{0,\alpha}$. 

In \cite{landais} it is proved that if $h \in C^{1,\alpha}(B_a)$. $f \in C^{0,\alpha}(B_a)$ and the distributional equation
\[ \int_{B_a} \frac{\nabla \psi\cdot \nabla h}{\sqrt{1+|\nabla h|^2}} = \int_{B_a} f\psi, \forall \psi \in C^1(B_a)\]
holds, then $h \in C^{2,\alpha}(B_a)$. Using this result we conclude that $h \in C^{2,\alpha}$ and $h$ is a strong solution of the equation
\begin{equation} -\di \left(\frac{\nabla h}{\sqrt{1+|\nabla h|^2}} \right) = \sum_{i=1}^m \mu_i |\nabla \mathfrak{u}_i|^2
\label{weakoptrel}
\end{equation}
with $\mu_1+...+\mu_p = 1$. Using straightforward computations equation \eqref{weakoptrel} is equivalent to
\[-\frac{\Delta h}{\sqrt{1+|\nabla h|^2}}+\frac{\nabla h\cdot \nabla h}{(1+|\nabla h|^2)^{3/2}}=\sum_{i=1}^m \mu_i |\nabla \mathfrak{u}_i|^2\]
which leads to
\begin{equation} -\Delta h = \sqrt{1+|\nabla h|^2} \left(\sum_{i=1}^m \mu_i |\nabla \mathfrak{u}_i|^2\right)-\frac{|\nabla h|^2}{1+|\nabla h|^2}.
\label{strongeq}
\end{equation}

Since $h \in C^{2,\alpha}(B_a)$ and this is true for any local chart which does not intersect $\Sigma$, the set $\Omega$ is a domain with $C^{2,\alpha}$ boundary outside $\Sigma$. Therefore using again the standard Schauder regularity results \cite[Theorem 6.18]{gilbarg-trudinger} we see that the eigenfunctions $(\mathfrak u_i)_{i=1}^p$ are in $C^{2,\alpha}(B_a)$. This means that the right hand side of equation \eqref{strongeq} is in $C^{1,\alpha}(B_a)$. Theorem 9.19 from \cite{gilbarg-trudinger} allows us to deduce that $h$ belongs to $C^{3,\alpha}(B_a)$. In general we see that if $h \in C^{k,\alpha}(B_a)$, $k \geq 2$, then the right hand side of \eqref{strongeq} is in $C^{k-1,\alpha}(B_a)$ and thus $h \in C^{k+1,\alpha}(B_a)$. An inductive bootstrap argument allows us to see that $h \in C^\infty(B_a)$. Moreover, since the coefficients of the partial differential equation \eqref{strongeq} are analytic and $h$ is $C^\infty$ the results of Morrey \cite{morrey-analytic} allow us to deduce that $h$ is analytic. Repeating the argument around each point outside $\Sigma$, we conclude that $\partial\Omega\setminus \Sigma$ is analytic. \hfill $\square$

\begin{rem}
\normalfont
Theorem \ref{main-regularity} improves the results in \cite{boouper} in yet another direction. In the cited reference the authors prove that at the optimum there exists a family of eigenfunctions $(\varphi_i)_{i=1}^m$ in the eigenspace associated to $\lambda_k$ such that
\[ \mathcal{H} = \sum_{i=1}^m (\partial_n \varphi_i)^2.\]
In the proof of our result we deduce a more precise result which says more about the family $(\varphi_i)$ and about the number of eigenfunctions present in the optimality condition. Indeed we prove that if $(\mathfrak{u}_i)_{i=1}^m$ is an orthonormal basis of the eigenspace associated to $\lambda_k$ then there exist $\mu_1,...,\mu_m \in [0,1]$ with $\mu_1+...+\mu_m=1$ such that
\[ \mathcal{H} = \sum_{i=1}^m \mu_i(\partial_n \mathfrak u_i)^2.\]
We note that this behavior has been conjectured by the numerical results presented in \cite{buboou}.
\end{rem}

\begin{rem}
\normalfont
The results of Theorem \ref{main-regularity} can be generalized to problems of the type
\[ \min_{\Omega \in \Bbb{R}^d} \left[F(\lambda_{k_1}(\Omega),...,\lambda_{k_p}(\Omega))+\Per(\Omega)\right] \]
where $F: \Bbb{R}^p \to \Bbb{R}_+$ satisfies the properties
\begin{itemize}
\item[(P1)] $F(x) \to \infty$ as $|x| \to \infty$.
\item[(P2)] $F$ is of class $C^1$ and at least one of its partial derivatives does not vanish when evaluated at $(\lambda_{k_1}(\Omega^*),...,\lambda_{k_p}(\Omega^*))$.
\item[(P3)] $F$ is increasing in each variable, furthermore, for any compact $K \subset \Bbb{R}^p\setminus \{0\}$ there exists $a>0$ such that if $x,y \in \Bbb{R}^p$ with $x_j \geq y_j,\ j=1,...,p$ then $F(x)-F(y) \geq a|x-y|$.
\end{itemize}
Properties (P1) and (P3) are taken from \cite{deveper} to guarantee existence. Property (P2) is stronger than the analogous property in \cite{deveper} (locally Lipschitz continuity) and allows us to differentiate the functional at the optimum. 
\end{rem}

\begin{rem}
\normalfont
The techniques presented here should also apply to the problem
\[ \min\{\lambda_k(\Omega)+\Per(\Omega), |\Omega| = m\}\]
introduced in \cite{ldpv}. Indeed, the functional is the same, but there is a volume constraint. The existence of a solution which is $C^{1,\alpha}$ regular outside a closed singular set of dimension at most $d-8$ has been treated in \cite{ldpv}. Like above it is possible to write the expression of the shape derivative, but we need to work with flows which preserve the volume. Such flows can be constructed for any perturbation $V$ such that $\int_{\partial \Omega} V.n = 0$ (see \cite{iliasoptim}). Then, using techniques similar to the ones in \cite{iliasoptim} we may deduce that at the optimality condition \eqref{weakoptrel} becomes
\begin{equation*} -\di \left(\frac{\nabla h}{\sqrt{1+|\nabla h|^2}} \right) = \sum_{i=1}^m \mu_i |\nabla \mathfrak{u}_i|^2+\text{const}.
\end{equation*}
This allows us to deduce the full regularity of the above problem in the same way as in Theorem \ref{main-regularity}.
\end{rem}

The results in this section allow us to give a different argument for \cite[Theorem 2.5]{bucbuthenper} and extend the corresponding result to all dimensions. The following result also generalises \cite[Theorem 5.7]{boouper}.

\begin{cor}
If $\Omega^*$ is a solution for problem \eqref{persum} then
\begin{itemize}
\item $\partial \Omega^*$ does not contain flat parts
\item $\partial \Omega^*$ does not contain open subsets of spheres unless $\Omega^*$ is a ball.
\end{itemize}
\end{cor}

\emph{Proof:} This is a simple consequence of the fact that $\partial \Omega^*$ is analytic outside a singular set of dimension at most $d-8$.

\section{The Multiplicity Cluster}

Once the regularity of the minimizers of \eqref{persum} is proved we may look into other facts about the optimal set. As noted in \cite{antunesf-per} and \cite{boouper} the numerical experiments show that in most cases the multiplicity of the $k$-th eigenvalue of the optimal set is larger than $1$. Numerical experiments also suggest that when $\lambda_k(\Omega^*)$ is multiple then $\lambda_k(\Omega^*)<\lambda_{k+1}(\Omega^*)$. A first result in this sense has been obtained in \cite{boouper} and is similar to \cite[Lemma 2.5.9]{henroteigs}. For the sake of completness we give a sketch of the proof.

\begin{thm}
Suppose $\Omega^*$ is a minimizer for problem \eqref{persum}. If we have $\lambda_{k-1}(\Omega^*)<\lambda_k(\Omega^*)$ then the $k$-th eigenvalue is simple, i.e. $\lambda_k(\Omega^*)<\lambda_{k+1}(\Omega^*)$.
\label{mcluster}
\end{thm}

\emph{Proof:} Suppose $\lambda_k(\Omega^*)$ is multiple and $\lambda_k(\Omega^*) = ... = \lambda_{k+m-1}(\Omega^*)$. Since $\lambda_k$ is the smallest eigenvalue in the cluster and $\Omega^*$ is a local minimizer, the $m$ analytic parametrizations of the eigenvalues of the spectrum, described in Section \ref{perturbations}, have derivative zero at $\varepsilon=0$. This would imply that $(\partial_n \mathfrak u_i)^2 = \mathcal H$ for all $i = 1,...,m$.

If $m>1$ this implies the existence of two different eigenfunctions $u,v$ associated to $\lambda_k(\Omega^*)$ such that $(\partial_n u)^2=(\partial_n v)^2 = \mathcal{H}$. The H\"olmgren uniqueness theorem implies that $\mathcal{H}$ cannot vanish on a relatively open set of $\partial \Omega$ (see \cite[Theorem 5.7]{boouper} for an alternative argument). Thus, there exists a relatively open set $\gamma \subset \partial \Omega$ such that $\mathcal{H}>0$ on $\gamma$. We can see that in this situation $\partial_n u = \pm \sqrt{\mathcal{H}}$ and $\partial_n v = \pm \sqrt{\mathcal{H}}$ on $\gamma$. This implies that we have either $\partial_n (u+v) = 0$ or $\partial_n (u-v) = 0$ on $\gamma$. By H\"olmgren's uniqueness theorem we deduce that $u=\pm v$ which is a contradiction. Therefore $\lambda_k(\Omega^*)$ is simple. \hfill $\square$

We can provide a stronger result based on the methods used in Theorem 2.5.8 and Lemma 2.5.9 from \cite{henroteigs}. We note that we have already proved that optimal sets $\Omega^*$ are analytic outside a singular set $\Sigma$ so the cited results apply in our case. Below we prove a generalization of Lemma 2.5.9 from \cite{henroteigs}.

\begin{lemma}
Suppose $\Omega$ is a smooth set in $\Bbb{R}^d$ and that $\lambda_{k+1}(\Omega)=...=\lambda_{k+m}(\Omega)$. Then there exist vector fields $V$ such that $\ell$ of  the derivatives $\displaystyle \frac{d\lambda_{k+i}(\Omega_\varepsilon)}{d\varepsilon}\Bigg|_{\varepsilon=0}$ are strictly negative and the other $m-\ell$ are strictly positive, for $\ell \in \{1,2,...,m-1\}$. Moreover, $V$ can be chosen in such a way that the derivative of the perimeter of $\Omega$ in the direction of $V$ is zero, which means that $\ds \int_{\partial \Omega} \mathcal{H} V.n d\sigma = 0$.
\label{signature}
\end{lemma}

\begin{rem}
Working like in \cite[Lemma 2.6]{boouper} we may construct a family of sets $\Omega_\varepsilon=\phi_\varepsilon(\Omega)$ such that $\Omega_\varepsilon$ have the same perimeter as $\Omega = \Omega_0$ and the derivative of the flow $\phi_\varepsilon$ at $\varepsilon=0$ is $V$. 
\label{fixed_per}
\end{rem}
\emph{Proof of Lemma \ref{signature}:} Theorem 2.5.8 from \cite{henroteigs} says that the derivatives of the eigenvalues of the perturbations $\Omega_\varepsilon$ of $\Omega$ are the eigenvalues of the $m\times m$ matrix
\[ \mathcal{A}_V = (a_{i,j}) \text{ where } a_{i,j} = -\int_{\partial \Omega} \partial_n \mathfrak u_i \partial_n \mathfrak u_j V.n d\sigma,\ \ 1 \leq i,j \leq m.\]
As a consequence, the derivatives of $\lambda_k(\Omega^*)+\Per(\Omega^*)$ in the direction $V$ are the eigenvalues of $\mathcal{A}$ shifted by the derivative of the perimeter for the same direction $V$.

As in the proof of the previous theorem we place ourselves in a relatively open subset $\gamma$ of $\partial \Omega$ where the mean curvature $\mathcal{H}$ is strictly positive. For a point $X \in \gamma$ we denote $\psi(X) = (\partial_n \mathfrak u_1(X),...,\partial_n \mathfrak u_m(X))\in \Bbb{R}^m$. We can find points $X_1,...,X_n \in \gamma$ such that the vectors $\psi(X_i),\ i=1,...,m$ are linearly independent. We present a simple proof of this fact in Appendix \ref{indep}.

Once chosen $X_1,...,X_n$ we consider a vector field $V_\delta$ such that $V_\delta$ is supported on $O_1\cup... \cup O_n$, where $X_p \in O_p$ and $\mathcal{H}^{n-1}(O_p) = \delta$. Given $j\in \{1,2,...,n\}$ we can consider $V_\delta$ to be of the form
\[ V_\delta = \begin{cases} (k_p/\delta)\vec n(x)  & \text{ on } O_p \setminus Q_p \\
 (f_p(x)/\delta)\vec n(x) & \text{ on }Q_p
\end{cases}\]
where constants $k_p$ satisfy $k_p < 0$ for $p\leq \ell$ and $k_p>0$ for $p > \ell$, $Q_p$ are portions around the boundary of $O_p$ of $\mathcal{H}^{d-1}$ measure $O(\delta^2)$, $f_p$ smooth extensions of the constant such that $f_p \leq k_p$ and $\vec n(x)$ denotes the normal to $\partial \Omega$. Furthermore, the constants $k_i$ are chosen such that $\ds \int_{\partial \Omega} \mathcal{H} V_\delta .nd\sigma = 0$, so that the derivative of the perimeter is zero. This can always be done. It suffices to fix the first $\ell$ constants $k_i = c<0$ and choose the rest equal to a positive constant which makes the above integral zero. Note that the derivative of the perimeter cannot be equal to zero if we choose all $k_i$ to have the same sign. This is due to the fact that $\gamma$ is a region where $\mathcal{H}$ is strictly positive. 

 It is not difficult to see that 
\[ \lim_{\delta \to 0} \int_{O_p} \partial_n \mathfrak u_i \partial_n \mathfrak u_j V_\delta.n d\sigma=k_p \partial_n \mathfrak u_i(X_i) \partial_n \mathfrak u_j(X_i), \ \ p = 1,...,n.\]
Therefore, the matrix $A_{V_\delta}$ converges to 
\[ k_1 \psi(X_i)^T\psi(X_i)+...+k_n \psi(X_n)^T \psi(X_n)\]
as $\delta$ goes to $0$.
For a general vector $X \in \Bbb{R}^n$ we have 
\[ X^T \mathcal A_{V_\delta} X \to k_1 (X.\psi(X_1))^2+...+k_n(X.\psi(X_n))^2. \]
Since the vectors $\psi(X_i)$ are linearly independent we deduce that this limit matrix  is associated to a quadratic form of signature $(n-\ell,\ell)$ so it has $\ell$ negative eigenvalues and $n-\ell$ positive ones. The continuity of the eigenvalues of the matrices with respect to the entries implies that for $\delta$ small enough $\mathcal{A}_{V_\delta}$ has $\ell$ negative eigenvalues and $n-\ell$ positive ones. Thus $\ell$ eigenvalues from the cluster have negative derivatives and $n-\ell$ have positive derivatives. \hfill $\square$

Note that even if the above result is stated for smooth sets, we may apply it in our case by performing only perturbations supported outside the singularity. An immediate corollary is the following.

\begin{cor}
If $\Omega^*$ is an optimizer for \eqref{persum} corresponding to the index $k$ then $\lambda_k(\Omega^*)<\lambda_{k+1}(\Omega^*)$. 
\label{last-cluster}
\end{cor}

\emph{Proof:} Suppose that $\lambda_k(\Omega^*)$ is multiple and $\lambda_s(\Omega^*)=...=\lambda_k(\Omega^*)=\lambda_{k+1}(\Omega^*)=...=\lambda_S(\Omega^*)$. Apply Lemma \ref{signature} to deduce that it is possible to find a perturbation $V$ of $\Omega^*$ such that the derivatives of the eigenvalues in the cluster are strictly negative for indices in $[s,k]$ and strictly positive for indices in $[k+1,S]$. Moreover, the perturbation can be chosen so that the perimeter has derivative zero in the direction of $V$. This perturbation would then decrease $\lambda_k+\Per$, thus contradicting optimality. Therefore we can only have $\lambda_k(\Omega^*)<\lambda_{k+1}(\Omega^*)$. \hfill $\square$

\begin{rem}
\normalfont
A direct consequence of the previous theorem is that the optimal cost for problem \eqref{persum} is a strictly increasing function of $k$, i.e.
\[ \min_{\Omega \in \Bbb{R}^d} \lambda_k(\Omega)+\Per(\Omega) < \min_{\Omega \in \Bbb{R}^d} \lambda_{k+1}(\Omega)+\Per(\Omega).\]
Denote $\Omega_k$ and $\Omega_{k+1}$ solutions of problem \eqref{persum} for $k$ and $k+1$, respectively. Then we have the following inequalities
\[ \lambda_k(\Omega_k)+\Per(\Omega_k) \leq \lambda_k(\Omega_{k+1})+\Per(\Omega_{k+1}) \leq \lambda_{k+1}(\Omega_{k+1})+\Per(\Omega_{k+1}),\]
where the first inequality comes from the optimality of $\Omega_k$ and the second from the ordering of the eigenvalues. If the optimal costs $c_k,\ c_{k+1}$ defined in \eqref{costk} satisfy $c_k=c_{k+1}$ then the above inequalities imply that $\lambda_k(\Omega_{k+1}) = \lambda_{k+1}(\Omega_{k+1})$ and that $\Omega_{k+1}$ is also solution of \eqref{persum} for $k$. This is in contradiction with the previous corollary.
\end{rem}

\begin{rem}
Remark \ref{fixed_per} allows us to conclude a similar result for the constrained problem
\[ d_k = \min \{\lambda_k(\Omega) : \Per(\Omega) = c\}. \]
In this case we also have $d_k<d_{k+1}$.
\end{rem}

\begin{rem}
\normalfont
Using the results above it is possible to deduce that for the shape $\Omega^*$ which minimizes $\lambda_2(\Omega)+\Per(\Omega)$ the second eigenvalue is simple and thus we have $(\partial_n u_2)^2 = \mathcal{H}$. Indeed, the optimal set $\Omega^*$ is connected \cite{deveper} so we have $\lambda_1(\Omega^*)<\lambda_2(\Omega^*)$. The results above show that $\lambda_2(\Omega^*)<\lambda_3(\Omega^*)$ so $\lambda_2(\Omega^*)$ is simple.
\end{rem}

\appendix

\section{From discrete to continous linear dependence}
\label{indep}

Consider a positive integer $n \geq 2$ and let $f_1,f_2,...,f_n$ be real continuous functions defined on an open connected set $I \subset \Bbb{R}^d$. For a point $x \in I$ we denote $\psi(x) = (f_1(x),...,f_n(x))\in \Bbb{R}^n$. We suppose that none of the functions $f_i$ is identically zero on $I$. Suppose that for each set of $n$ different points $x_1,...,x_n \in I$ 
the vectors
\[ \psi(x_1),...,\psi(x_n)\]
are linearly dependent on an open connected subset of $I$. Then the functions $f_1,...,f_n$ are linearly dependent, i.e. there exist constants $\alpha_1,...,\alpha_n$ not all zero such that 
\[ \alpha_1 f_1+...+\alpha_n f_n = 0.\]

Let's start with the case $n=2$. Pick $x_1 \neq x_2$ in $I$. We know that the vectors $\psi(x_1),\psi(x_2)$ are linearly dependent which means that 
\[ \begin{pmatrix}
 f_1(x_1) \\ f_2(x_1) 
\end{pmatrix}
=\lambda
\begin{pmatrix}
f_1(x_2) \\ f_2(x_2)
\end{pmatrix}
\]
We find that on the set $\{f_2 \neq 0\}$ the function $f_1/f_2$ is constant. Thus, picking eventually a connected component of $\{f_2 \neq 0\}$ we find that $f_1,f_2$ are linearly dependent on a connected open set.

We can prove the result for general $n \geq 2$ by induction. Indeed, let's suppose that the result holds for $n$ functions. Consider now $f_1,...,f_{n+1}$ defined on $I$ such that for every $n+1$ points $x_1,...,x_{n+1}$ the vectors $\psi(x_1),...,\psi(x_{n+1})$ are linearly dependent. Suppose that no two functions $f_i$ are equal or else we have nothing to prove. Furthermore, we choose $x_1$ such that $f(x_1) \neq 0$. Thus there exist scalars $\lambda_1,...,\lambda_n$, not all zero such that
\[ \begin{cases}
\hfill \lambda_1f_1(x_1) +\lambda_2 f_1(x_2)+...+\lambda_{n+1}f_1(x_{n+1}) & = 0 \\
\hfill \lambda_1f_2(x_1) +\lambda_2f_2(x_2)+...+\lambda_{n+1}f_2(x_{n+1}) & = 0 \\
\hfill \cdots &  \\
\lambda_1f_{n+1}(x_1) +\lambda_2 f_{n+1}(x_2)+...+\lambda_{n+1}f_{n+1}(x_{n+1}) & = 0. 
\end{cases}\]
Without loss of generality suppose that $\lambda_1 \neq 0$. Eliminating the elements below $\lambda_1f_1(x_1)$ in the first column leaves us with a system of the form
\[ 
\begin{cases}
\mu_2(f_2-f_1)(x_2)+...+\mu_{n+1}(f_{2}-f_1)(x_{n+1}) & = 0 \\
\cdots & = 0 \\
\mu_2(f_{n+1}-f_1)(x_2)+...+\mu_{n+1}(f_{n+1}-f_1)(x_{n+1}) & = 0 \\
\end{cases}
\]
Therefore $f_2-f_1,...,f_{n+1}-f_1$ satisfy the induction hypothesis and thus they are linearly dependent. It is obvious that this implies that $f_1,...,f_{n+1}$ are linearly dependent.

\section{A variant of the Hahn-Banach Theorem}
\label{hb}

We recall below the well known Hahn-Banach separation theorem. For a proof or more details see \cite[Section 2.4]{clarke}.
\begin{thm}
Ket $K_1,K_2$ be nonempty, disjoint convex subsets of the normed vector space $X$. Then if $K_1$ is open, there exist $\zeta \in X^*$ and $\theta \in \Bbb{R}$ such that
\[ \langle \zeta,x\rangle < \theta \leq \langle \zeta, y\rangle \forall x \in K_1,\ y \in K_2.\]
\label{thb}
\end{thm}
An immediate consequence is the following.
\begin{prop}
Let $\{\zeta_i : i=1,2,...,k\}$ be a finite subset in $X^*$. The following are equivalent:
\begin{itemize}
\item[(a)] There is no $\nu \in X$ such that $\langle \zeta_i,\nu \rangle < 0,\ \forall i \in \{1,2,...,k\}$.
\item[(b)] The set $\{\zeta_i : i=1,2,...,k\}$ is positively linearly dependent: there exists a non zero nonnegative vector $\gamma \in \Bbb{R}^k$ such that $\sum_{i=1}^k \gamma_i \zeta_i = 0$.
\end{itemize}
\label{hb-variant}
\end{prop}

\emph{Proof:} First let's note that we may suppose that none of the functionals $\zeta_i$ are identically zero, since the problem would reduce itself to a smaller value of $k$ or to the trivial case when all functionals considered vanish. 

We start with the implication $(b) \Longrightarrow (a)$. Suppose that $\sum_{i=1}^k \gamma_i \zeta_i = 0$. Now, if $\nu \in X$ we have
\[ 0 = \sum_{i=1}^k \gamma_i \langle \zeta_i,\nu \rangle \leq (\sum_{i=1}^k\gamma_i)\max_i \langle \zeta_i,\nu \rangle, \]
which implies $\max_i \langle \zeta_i,\nu \rangle \geq 0$, since $\sum_{i=1}^k \gamma_i>0$.

For the implication $(a) \Longrightarrow (b)$ we consider the following subsets of $\Bbb{R}^k$:
\[ K_1 = \{y \in \Bbb{R}^k : y_i<0 \forall i \in \{1,2,...,k\}\}\]
\[ K_2 = \{ \left( \langle \zeta_1,\nu \rangle, \langle \zeta_2,\nu \rangle,...,\langle \zeta_k,\nu \rangle\right) : \nu \in X\}.\]
We see immediately that (a) implies $K_1 \cap K_2 = \emptyset$ and furthermore $K_1$ is open. Therefore, we can apply the Hahn-Banach separation theorem \ref{thb} and find a functional $\varphi:\Bbb{R}^k \to \Bbb{R}$ and $\theta \in \Bbb{R}$ such that 
\[ \langle \varphi, x\rangle < \theta \leq \langle \varphi,y \rangle, \forall x \in K_1,y\in K_2.\]
We know that such a functional is of the form
\[ \varphi(x) = \sum_{i=1}^k \gamma_ix_i,\]
where $\gamma=(\gamma_1,...,\gamma_k)$ is a non-zero element of $\Bbb{R}^k$. Therefore $\langle \varphi , y \rangle \geq \theta, \forall y \in K_2$ becomes
\[ \theta \leq \langle \sum_{i=1}^k \gamma_i\zeta_i ,\nu \rangle, \forall \nu \in X,\]
which is only possible when $\sum_{i=1}^k \gamma_i \zeta_i=0$. On the other hand, we have 
\[ \sum_{i=1}^k \gamma_i y_i <\theta \text{ for all }y=(y_i) \in K_1.\]
Taking $y_i^n = -n,\ y_j^n=-1/n,\ j \neq i$ and taking $n \to \infty$ we see that this is possible only if all $\gamma_i$ are non-negative. This finishes the proof. \hfill $\square$

\section{A formula relating the Laplacian and the Hessian}
\label{int-parts}

Theorem 3.1.1.1. from \cite{grisvard} says that for $v \in H^1(\Omega)^n$ and $\Omega$ of class $C^2$ we have
\[ \int_\Omega |\di v|^2 dx -\sum_{i,j=1}^n \int_\Omega \partial_i v_j \partial_j v_i dx = \text{(tangential part) }+\int_{\partial \Omega} \mathcal{H} (v\cdot n\nu)^2 d\sigma, \]
where $\mathcal{H}$ is the mean curvature of $\Omega$. If we put $v = \nabla u$ for $u=0$ on $\partial \Omega$, then the tangential part is zero and we obtain a relation between the Laplacian and the Hessian. The aim here is to write a similar relation when we have an additional multiplication with a function $\varphi \in {H}^1(\Bbb{R}^d)$. 

We wish to prove the following
\begin{prop}
If $\Omega$ is of class $C^2$, $\varphi \in H^1(\Bbb{R}^d)$ and $u$ is smooth in $\Omega \cap \operatorname{supp} \varphi$ then 
\[ \int_\Omega \varphi |D^2u|^2dx + \int_{\partial \Omega} \varphi \mathcal{H}|\nabla u|^2d\sigma = \int_\Omega \varphi (\Delta u)^2 dx+\int_\Omega (\nabla u \cdot \nabla \varphi)\Delta u dx -\int_\Omega \nabla \varphi \cdot D^2uDudx \]
\label{int-by-parts}
\end{prop}

\emph{Proof:} We integrate by parts two times and we obtain
\begin{align*}
\int_\Omega \varphi |\di v|^2 & = \sum_{i,j=1}^n \int_\Omega \varphi \partial_i v_i \partial_j v_j dx \\
& \stackrel{\partial_i}{=} -\sum_{i,j=1}^n  \int_\Omega v_i \varphi \partial_{i,j} v_jdx -\sum_{i,j=1}^n \int_\Omega v_i \partial_i \varphi \partial_j v_jdx + \sum_{i,j=1}^n\int_{\partial \Omega} \varphi v_i \partial_j v_j \nu_i d\sigma  \\
& \stackrel{\partial_j}{=} \sum_{i,j=1}^n \int_\Omega \varphi \partial_j v_i \partial_i v_j dx +\sum_{i,j=1}^n\int_\Omega \partial_j \varphi v_i \partial_i v_j dx -\sum_{i,j=1}^n \int_{\partial \Omega} \varphi v_i \partial_i v_j \nu_jd\sigma  - \\
& \ \ \  -\sum_{i,j=1}^n \int_\Omega v_i \partial_i \varphi \partial_j v_jdx + \sum_{i,j=1}^n\int_{\partial \Omega} \varphi v_i \partial_j v_j \nu_i d\sigma 
\end{align*}
Therefore we have
\begin{align*}
I(v) & = \int_\Omega \varphi |\di v|^2 -\sum_{i,j=1}^n \int_\Omega \varphi \partial_j v_i \partial_i v_j dx\\
& = A_\Omega \text{ (terms on }\Omega) + A_{\partial \Omega} \text{ (terms on }\partial \Omega)
\end{align*}

The terms integrated on $\Omega$ are
\begin{align*}
A_\Omega & = \sum_{i,j=1}^n\int_\Omega \partial_j \varphi v_i \partial_i v_j dx  -\sum_{i,j=1}^n \int_\Omega v_i \partial_i \varphi \partial_j v_jdx \\
& = \int_\Omega \{ (v \cdot \nabla)v \cdot \nabla \varphi\}dx - \int_\Omega (v\cdot \nabla \varphi)\di v dx.
\end{align*}

In order to avoid possible complications, unnecessary in our case, we recall that we want to use this computation for $v = \nabla u$ with $u=0$ on $\partial \Omega$. This means, in particular, that the tangential part of $v$ is $v_\tau = (\nabla u)_\tau = 0$. Therefore, in the following, we suppose that the tangential part of $v$ is zero. Using the reasoning from \cite[Theorem 3.1.1.1]{grisvard} and neglecting the tangential components we obtain
\begin{align*}A_{\partial \Omega} & =   -\sum_{i,j=1}^n \int_{\partial \Omega} \varphi v_i \partial_i v_j \nu_j d\sigma+ \sum_{i,j=1}^n\int_{\partial \Omega} \varphi v_i \partial_j v_j \nu_i d\sigma \\
& = \int_{\partial \Omega} \varphi v_\nu \di v d\sigma - \int_{\partial \Omega} \varphi\{(v\cdot \nabla)v\}\cdot \nu d\sigma \\
& = \int_{\partial \Omega} \varphi \mathcal{H} v^2d\sigma 
\end{align*}

As stated before, we replace $v$ by $\nabla u$ where $u=0$ on $\partial \Omega$ and we get
\[ \int_\Omega \varphi(\Delta u)^2dx -\int_\Omega \varphi |D^2u|^2dx = \int_\Omega \nabla \varphi \cdot D^2uDu dx- \int_\Omega (\nabla u \cdot \nabla \varphi)\Delta u dx+ \int_{\partial \Omega} \varphi \mathcal{H} |\nabla u|^2d\sigma.  \]
Therefore
\[ \int_\Omega \varphi |D^2u|^2dx + \int_{\partial \Omega} \varphi \mathcal{H}|\nabla u|^2d\sigma = \int_\Omega \varphi (\Delta u)^2 dx+\int_\Omega (\nabla u \cdot \nabla \varphi)\Delta u dx -\int_\Omega \nabla \varphi \cdot D^2uDudx \] \hfill $\square$

Note that even if the result above is stated for $C^2$ sets, we may apply it in the case where $\Omega$ is $C^{1,\alpha}$ with non-negative and bounded distributional curvature. Analysing the proof we notice that the only difficulty comes from the boundary integrals $A_{\partial \Omega}$ where which are considered for fields $v$ with zero tangential components.

\begin{ack}
The author wishes to thank Dorin Bucur, Otis Codosh, Jimmy \mbox{Lamboley} and Bozhidar Velichkov for useful discussions in connection with this work. The author also thanks the anonymous reviewers for their careful analysis of the text and for the suggestions which helped improve the quality of this paper. This work was supported by the project ANR Optiform and by the Fondation Sciences Mathematiques de Paris.
\end{ack}
\bibliographystyle{abbrv}
\bibliography{../master.bib}

\Addresses
\end{document}